\documentclass[3p,12pt]{elsarticle}

\usepackage{amssymb}

\journal{Stochastic Analysis and Applications}
\newtheorem{thm}{Theorem}
\newdefinition{rmk}{Remark}

\begin{document}

\begin{frontmatter}

\title{\textbf{Uniform convergence of wavelet expansions}\\ \textbf{of Gaussian random processes}\\[5mm]
\small \textit{Short title: Uniform convergence of wavelet expansions}}

\author[ku]{Yuriy Kozachenko}
\ead{ykoz@ukr.net}

\author[la]{Andriy Olenko\corref{cor1}}
\ead{a.olenko@latrobe.edu.au}

\author[ku]{Olga Polosmak}
\ead{DidenkoOlga@yandex.ru}

\cortext[cor1]{Corresponding author. Phone: +61-3-9479-2609 \quad  Fax:  +61-3-9479-2466}
\address[ku]{Department of Probability Theory, Statistics and Actuarial Mathematics, Kyiv University, Kyiv, Ukraine}
\address[la]{Department of Mathematics and Statistics, La Trobe University, Victoria 3086, Australia\\

\

 This is an Author's Accepted Manuscript of an article published in the
Stochastic Analysis and Applications, Vol. 29, No. 2, 169--184.  [copyright Taylor \& Francis], available online at: http://www.tandfonline.com/ [DOI:10.1080/07362994.2011.532034]
}

\begin{abstract}
New results on uniform
convergence in probability for the most general classes of wavelet expansions of
stationary Gaussian random processes are given.
\end{abstract}


\end{frontmatter}


\section{Introduction}
In various applications in data compression, signal processing and simulation, it could be useful to convert the problem of analyzing a continuous-time random process to that of analyzing a random sequence, which is much simpler. Multi\-resolution analysis provides an efficient framework for the decomposition of random processes. This approach is widely used in statistics to estimate a curve given observations of the curve plus some noise.

Various extensions of the standard statistical methodology were proposed recently. These include curve estimation in the presence of correlated noise. For these purposes the wavelet based expansions have numerous advantages over Fourier series, see \cite{kur, wal}, and often lead to stable computations, see \cite{pho}.

However, in many cases numerical simulation results need to be confirmed by theoretical analysis. Recently, a considerable attention was given
to the properties of the wavelet transform and of the wavelet orthonormal series representation of
random processes. More information on convergence of wavelet expansions of random processes in various spaces, references and numerous applications can be found in \cite{cam, did, ist, kozvas, kozpol, kur, zha}.

We focus out attention on  uniform convergence of wavelet
expansions for stationary Gaussian random processes. We consider a random process $\mathbf{X}(t)$ such that $\mathbf E\mathbf{X}(t)=0$ for all $t\in\mathbb R$.

In the book \cite{har} wavelet expansions of functions bounded on $\mathbb R$ were studied in different
spaces. Obtained results were applied by several authors to investigate wavelet
expansions of random processes $\mathbf{X}(t)$ bounded on $\mathbb R$. However
the majority of random processes, which are interesting from theoretical and practical application points of view,
has almost surely unbounded sample paths on $\mathbb R$. In numerous cases developed deterministic methods may not be appropriate to investigate wavelet expansions of stochastic processes.  It indicates the necessity of elaborating special stochastic techniques.

In the paper we consider  stationary Gaussian random processes $\mathbf{X}(t)$ and their approximations by sums of wavelet functions
\begin{equation}\label{Xn}\mathbf{X}_{n,\mathbf{k}_n}(t):=\sum_{|k|\le k_0}\xi_{0k}\phi_{0k}(t)+\sum_{j=0}^{n-1}\sum_{|k|\le k_j}\eta_{jk}\psi_{jk}(t)\,,
\end{equation}
where $\mathbf{k}_n:=(k_0,...,k_{n-1}).$

Contrary to many theoretical results (see, for example, \cite{kozvas,kur}) with infinite series form of $\mathbf{X}_{n,\mathbf{k}_n}(t),$ in direct numerical implementations we always consider truncated series like (\ref{Xn}), where the number of terms in the sums is finite by application reasons. However, there are almost no stochastic results on uniform convergence of finite wavelet expansions to $\mathbf{X}(t).$

We show that, under suitable conditions, the sequence $\mathbf{X}_{n,\mathbf{k}_n}(t)$ converges in probability in Banach space $C([0,T])$, i.e.
$$P\left\{\sup_{0\le t\le T} |\mathbf{X}(t)-\mathbf{X}_{n,\mathbf{k}_n}(t)|>\varepsilon \right\}\to 0,
$$
when  $n\to\infty$ and $k_j\to\infty$ for all $j\in \mathbb{N}_0:=\{0,1,...\}\,.$ More details on the general theory of random processes in the space $C(\mathbb T)$ can be found in \cite{bulkoz}.

The numbers $n$ and $\mathbf{k}_n$ of terms in the truncated series $\mathbf{X}_{n,\mathbf{k}_n}(t)$ can approach infinity in any arbitrary way. Thought about this way, one sees that the paper deals with the most general class of such wavelet expansions in comparison with particular cases considered by different authors, see, for example, \cite{cam,kur}.

Most known results (see, for example, \cite{cam,did,ist,wong,zha}) concern the mean-square convergence, but for practical applications one needs to require uniform convergence.

We present the first result on stochastic uniform convergence of general finite wavelet expansions in the open literature.
It should be mentioned that in the general case the random coefficients in (\ref{Xn}) may be dependent and form an overcomplete system of basis functions. That is why the implementation of the proposed method has promising potential for nonstationary random processes.

The organization of this article is the following. In the
second section we introduce the necessary background from wavelet
theory and certain sufficient conditions for mean-square convergence of wavelet expansions in the space $L_2(\Omega),$
obtained in~\cite{kozpol}. In \S 3 we formulate and discuss the main theorem on uniform convergence in probability of the wavelet expansions of stationary Gaussian random processes. The next section contains the proof of the theorem. Conclusions are made in section 5.

\section{Wavelet representation of random processes}
Let $\phi(x),$ $x\in\mathbb R$ be a function from the space
$L_2(\mathbb R)$ such that $\widehat{\phi}(0)\ne 0$ and  $\widehat{\phi}(y)$ is
continuous at $0,$ where
$$\widehat{\phi}(y)=\int_{\mathbb
R}e^{-iyx}{\phi(x)}\,dx$$
is the Fourier transform of $\phi.$

 Suppose that the following assumption
holds true:
$$\sum_{k\in\mathbb Z} |\widehat{\phi}(y+2{\pi}k)|^2=1\  {\rm (a.e.)}
$$
There exists a function $m_0(x)\in L_2([0,2\pi])$, such that $m_0(x)$
has the period $2\pi$ and
$$\widehat{\phi}(y)=m_0\left(y/2\right)\widehat{\phi}\left(y/2\right)\ {\rm (a.e.)}
$$
 In this case the  function $\phi(x)$ is called
the $f$-wavelet.

Let $\psi(x)$ be the inverse Fourier transform of the function
$$\widehat{\psi}(y)=\overline{m_0\left(\frac
y2+\pi\right)}\cdot\exp\left\{-i\frac
y2\right\}\cdot\widehat{\phi}\left(\frac y2\right).$$
Then the function
$$\psi(x)=\frac1{2\pi}\int_{\mathbb
R}e^{iyx}{\widehat{\psi}(y)}\,dy$$ is called the $m$-wavelet.

Let
\begin{equation}\label{2phijk}\phi_{jk}(x)=2^{j/2}\phi(2^jx-k),\quad
\psi_{jk}(x)=2^{j/2}\psi(2^jx-k),\quad j,k \in\mathbb Z\,.
\end{equation}
It is known that the family of functions $\{\phi_{0k};
\psi_{jk},j\in \mathbb N_0\}$ is an orthonormal basis in
$L_2(\mathbb R)$ (see, for example, \cite{chu,dau}).

An arbitrary function $f(x)\in L_2(\mathbb R)$ can be represented in the form
\begin{equation}\label{2.5}f(x)=\sum_{k\in\mathbb Z}\alpha_{0k}\phi_{0k}(x)+\sum_{j=0}^{\infty}\sum_{k\in\mathbb Z}\beta_{jk}\psi_{jk}(x)\,,
\end{equation}
$$\alpha_{0k}=\int_{\mathbb R}f(x)\overline{\phi_{0k}(x)}\,dx,\quad \beta_{jk}=\int_{\mathbb R}f(x)\overline{\psi_{jk}(x)}\,dx.$$
The representation (\ref{2.5}) is called a wavelet representation.

The series~(\ref{2.5}) converges in the space
$L_2(\mathbb R),$ i.e.
$\sum_{k\in\mathbb Z}|\alpha_{0k}|^2+\sum_{j=0}^{\infty}\sum_{k\in\mathbb Z}|\beta_{jk}|^2<\infty\,.$

The integrals $\alpha_{0k}$ and $\beta_{jk}$ may also exist
for functions from $L_1(\mathbb R)$ and other function spaces. Therefore it is possible to obtain the representation~(\ref{2.5}) for function classes which are wider than $L_2(\mathbb R)$.

Let $\{\Omega, \cal{B}, P\}$ be a standard probability space. Let
$\mathbf{X}(t),$ $t\in\mathbb R$ be a  random process
such that $\mathbf E\mathbf{X}(t)=0\,.$ It is possible to obtain representations like~(\ref{2.5}) for random
processes, if sample trajectories of these processes are  in the space $L_2(\mathbb R).$ However  the majority of random
processes do not possess this property. For example, sample paths of
stationary processes are not in the space $L_2(\mathbb R)$ (a.s.).

We want to construct a representation of the kind~(\ref{2.5}) for
$\mathbf{X}(t)$ with mean-square integrals
$$\xi_{0k}=\int_{\mathbb R}\mathbf{X}(t)\overline{\phi_{0k}(t)}\,dt,\quad \eta_{jk}=\int_{\mathbb R}\mathbf{X}(t)\overline{\psi_{jk}(t)}\,dt\,.$$
Consider the approximants $\mathbf{X}_{n,\mathbf{k}_n}(t)$ of $\mathbf{X}(t)$ defined by (\ref{Xn}). Theorem~\ref{8103} below guarantees the mean-square convergence of $\mathbf{X}_{n,\mathbf{k}_n}(t)$  to $\mathbf{X}(t).$\vspace{1mm}

\noindent {\bf Assumption S.} \cite{har}  For the $f$-wavelet $\phi$ there exists a function
$\Phi(x),$ $x\ge 0$ such that $\Phi(0)<\infty$, $\Phi(x)$
is a decreasing function, $|\phi(x)|\le \Phi(|x|)$ (a.e.) and $\int_{\mathbb R}\Phi(|x|)\,dx<\infty\,.$

Let $c(x),$ $x\in\mathbb R$ denote a non decreasing  even function on $[0,\infty)$ with $c(0)>0.$

\begin{thm}\label{8103} {\rm\cite{kozpol}} Let $\mathbf{X}(t),$ $t\in\mathbb R$ be a random process
such that $\mathbf E\mathbf{X}(t)=0,$ $\mathbf E|\mathbf{X}(t)|^2<\infty$ for all $t\in\mathbb R,$ and
its covariance function  $R(t,s)$  is continuous. Let the $f$-wavelet $\phi$
and the $m$-wavelet $\psi$ be continuous functions and  the assumption {\rm S}  hold true
for both  $\phi$ and  $\psi.$
Suppose that there exists a
function $A:(0,\infty)\to (0,\infty)$ and $x_0\in\mathbb R$ such that  $
c(ax)\le c(x)\cdot A(a),$ for all $x\ge x_0.$ 

If
$$\int\limits_{\mathbb R}c(x)\Phi(|x|)\,dx<\infty\quad {\rm and} \quad
|R(t,t)|^{1/2}\le c(t)$$
then
\begin{enumerate}
  \item[\rm 1.]  $\mathbf{X}_{n,\mathbf{k}_n}(t)\in L_2(\Omega)\,;$
  \item[\rm 2.]  $\mathbf{X}_{n,\mathbf{k}_n}(t)\to \mathbf{X}(t)$ in mean square when  $n\to\infty$ and $k_j\to\infty$ for all
$j\in \mathbb N_0\,.$
\end{enumerate}
\end{thm}
\begin{rmk}
For stationary Gaussian processes  we can choose $c^2(x)\equiv R(0,0)$ and~$A(a)\equiv 1.$
\end{rmk}

\section{Uniform convergence  of wavelet expansions
for Gaussian random processes}
\begin{thm}\label{213}{\rm \cite{kozsli}} Let $\mathbb T =[0,T]$, $\rho(t,s)=|t-s|$. Let $\mathbf{X}_n(t),$ $t\in\mathbb T$
 be a sequence of Gaussian stochastic processes. Assume that
all  $\mathbf{X}_n(t)$ are separable in $(\mathbb T, \rho)$ and
$$\sup_{n\ge 1}\sup_{|t-s|\le h}\left( \mathbf E|\mathbf{X}_n(t)-\mathbf{X}_n(s)|^2\right)^{1/2}\le \sigma(h)\,,$$
where $\sigma(h)$ is a monotone increasing function such that
$\sigma(h)\to 0$ when $h\to 0\,.$

Suppose that for some $\varepsilon>0$
\begin{equation}\label{142}\int\limits_0^\varepsilon
\sqrt{-\ln\left(\sigma^{(-1)}(u)\right)}\,du<\infty\,,\end{equation}
where $\sigma^{(-1)}(u)$ is the inverse function of $\sigma(u)$.
 If the processes $\mathbf{X}_n(t)$ converges in probability to the process
 $\mathbf{X}(t)$ for all $t\in \mathbb T$, then $\mathbf{X}_n(t)$ converges in
 probability to  $\mathbf{X}(t)$ in the space
$C(\mathbb T).$
\end{thm}
\begin{rmk}\label{772}For example, it is easy to check that the assumption~(\ref{142})
holds true for  $$\sigma(h)=\frac c{\left(\ln(e^{\alpha}+
\frac1{h})\right)^{\alpha}}\quad \mbox {and} \quad \sigma(h)=c h^{\gamma},$$
 when $c>0,$ $\alpha>1/2,$ $\gamma>0.$
\end{rmk}

Now we are ready to formulate the main result.
\begin{thm}\label{main} Let $\mathbf{X}(t),$ $t\in\mathbb R$ be a stationary separable centered Gaussian random process such that
 its covariance function  $R(t,s)=R(t-s)$ is continuous. Let the $f$-wavelet $\phi$
and the corresponding $m$-wavelet $\psi$ be continuous functions and  the assumption {\rm S}  hold true
for both  $\phi$ and  $\psi.$ Suppose that the following
conditions hold:

{\rm \begin{enumerate}
  \item\label{con1} there exist $\phi'(u),$ $\widehat{\psi}''(u),$ and  $\widehat{\psi}(0)=0,$ $\widehat{\psi}'(0)=0 ;$
  \item\label{con2} $c_{\phi}:=\sup\limits_{u\in \mathbb R}|\widehat{\phi}(u)|<\infty,$ $
c_{\phi'}:=\sup\limits_{u\in \mathbb R}|\widehat{\phi}'(u)|<\infty,$ $c_{\psi''}:=\sup\limits_{u\in \mathbb R}|\widehat{\psi}''(u)|<\infty;$

  \item $\widehat{\phi}(u)\to 0$ and $\widehat{\psi}(u)\to 0$ when $u\to \pm\infty;$
  \item\label{con4} there exist $0<\gamma<\frac{1}{2}$  and $\alpha>\frac{1}{2}$ such that $\int\limits_{\mathbb
R}\left(\ln(1+|u|)\right)^{\alpha}|\widehat{\psi}(u)|^{\gamma}\,du<\infty,$\\
 $\int\limits_{\mathbb R}\left(\ln(1+|u|)\right)^{\alpha}|\widehat{\phi}(u)|^{\gamma}\,du<\infty\,;$

  \item\label{con5} there exists $\widehat R(z)$ and  $\sup\limits_{z\in \mathbb R}\left|\widehat R(z)\right|<\infty\,;$

  \item\label{con6} $\int\limits_{\mathbb R}\left|\widehat
R'(z)\right|\,dz<\infty$ and $\int\limits_{\mathbb R}\left|\widehat
R^{(p)}(z)\right||z|^4\,dz<\infty$ for $p=0,1\,.$
\end{enumerate}}
\noindent Then $\mathbf{X}_{n,\mathbf{k}_n}(t)\to \mathbf{X}(t)$ uniformly in probability on each interval $[0,T]$ when  $n\to\infty$ and $k_j\to\infty$ for all
$j\in \mathbb N_0\,.$
\end{thm}

Before seeing the proof of the theorem, we clarify the role of some of the assumptions.  Two kinds of assumptions were made:
\begin{itemize}
  \item conditions \ref{con1}-\ref{con4} on the wavelet basis and
  \item conditions \ref{con5} and \ref{con6} on the random process.
\end{itemize}

Contrary to many other results in literature, our assumptions are very simple and can be easily verified.

Conditions \ref{con1}-\ref{con4} are related to the smoothness and the decay rate of the wavelet basis functions $\phi$ and $\psi.$ It is easy to check that numerous wavelets satisfy these conditions, for example, the well known Daubechies, Battle-Lemarie and Meyer wavelet bases. Conditions \ref{con5} and \ref{con6} on the random process $\mathbf{X}(t)$ are formulated in terms of the spectral density $\widehat R(z).$ These conditions are related to the behavior of the high-frequency part of the spectrum. Both sets of assumptions are standard in the convergence studies.

If we narrow our general class of wavelet expansions and impose some additional constraints on rates of the sequences $\mathbf{k}_n$ we can enlarge classes of wavelets bases and random processes in the theorem.  It will be seen from the proof of the theorem that all we need are conditions which guarantee that the series $Q,$ $B_1,$ $q_1,$ $q_{\phi},$ and $q_{\phi 1}$ are convergent.

\section{Proof of the main theorem}
From conditions \ref{con2} and \ref{con4}, it follows that
$$\int\limits_{\mathbb
R}\left(\ln(1+|u|)\right)^{\alpha_1}|\widehat{\psi}(u)|^{\gamma}\,du<\infty\quad \mbox{and}\quad \int\limits_{\mathbb
R}\left(\ln(1+|u|)\right)^{\alpha_1}|\widehat{\phi}(u)|^{\gamma}\,du<\infty$$
 for any $\alpha_1\in(1/2,\alpha).$ Thus to prove the
theorem we only  consider the case $1/2<\alpha\le 1\,.$

We first prove  that for some $B>0$ the inequality
\begin{equation}\label{1075}\left(\mathbf E\left|\mathbf{X}_{n,\mathbf{k}_n}(t)- \mathbf{X}_{n,\mathbf{k}_n}(s)\right|^2\right)^{1/2}\le \frac B{\left(\ln\left(e^{\alpha}+
\frac1{|t-s|}\right)\right)^{\alpha}}, {\  } 1/2<\alpha\le 1
\end{equation}
holds true for all $n\in \mathbb N,$ $\mathbf{k}_n\in \mathbb N_0^n,$ and $t,s\in [0,T].$

By (\ref{Xn}) we obtain
$$(\mathbf E\left|\mathbf{X}_{n,\mathbf{k}_n}(t)-
\mathbf{X}_{n,\mathbf{k}_n}(s)\right|^2)^{1/2}\le \left(\mathbf E\left|\sum\limits_{|k|\le
k_0}\xi_{0k}(\phi_{0k}(t)-\phi_{0k}(s))\right|^2\right)^{1/2}$$
$$+\sum\limits_{j=0}^{n-1}\left(\mathbf E\left|\sum\limits_{|k|\le
k_j}\eta_{jk}(\psi_{jk}(t)-\psi_{jk}(s))\right|^2\right)^{1/2}=:\sqrt{S}+\sum\limits_{j=0}^{n-1}\sqrt{S_j}\,.$$

We will show how to handle $S_j$, then similar techniques can be used
to deal with the remaining term $S.$

$S_j$ satisfies the inequality
 $$S_j\le \sum\limits_{|k|\le k_j}\sum\limits_{|l|\le k_j}|\mathbf E\eta_{jk}\overline{\eta_{jl}}||\psi_{jk}(t)-\psi_{jk}(s)||\psi_{jl}(t)-\psi_{jl}(s)|\,.$$
Let us consider $\mathbf E\eta_{jk}\overline{\eta_{jl}}.$ By means of Parseval's theorem
we deduce
$${\  }\mathbf E\eta_{jk}\overline{\eta_{jl}}=\int\limits_{\mathbb
R}\int\limits_{\mathbb
R}\mathbf E \mathbf{X}(u)\overline{\mathbf{X}(v)}\ \overline{\psi_{jk}(u)}\psi_{jl}(v)\,dudv=\int\limits_{\mathbb R}\int\limits_{\mathbb R}R(u-v)\overline{\psi_{jk}(u)}\,du\,\psi_{jl}(v)\,dv$$
$$=\int\limits_{\mathbb R}\frac{1}{2 \pi}\int\limits_{\mathbb R}\widehat{R}(z)e^{-ivz}\overline{\widehat{\psi}_{jk}(z)}\,dz\psi_{jl}(v)\,dv=\frac{1}{2 \pi}\int\limits_{\mathbb R}\widehat{R}(z)\,\overline{\widehat{\psi}_{jk}(z)}\,\widehat{\psi}_{jl}(z)\,dz\,.$$

The order of integration can be changed because
$$ \int\limits_{\mathbb R}\int\limits_{\mathbb R}\left|\widehat{R}(z)e^{-ivz}\overline{\widehat{\psi}_{jk}(z)}\psi_{jl}(v)\right|\,dz\,dv\le \sup_{z\in\mathbb R}|\widehat{R}(z)|\cdot
\int\limits_{\mathbb R}|{\widehat{\psi}_{jk}(z)}|\,dz\cdot\int\limits_{\mathbb R}|\psi_{jl}(v)|\,dv<\infty\,.$$
The last expression is finite due to  (\ref{2phijk}), the assumption S, the estimate (\ref{c2}), and the representation
\begin{equation}\label{104.11}\widehat{\psi}_{jk}(z)=\frac{e^{-i\frac k{2^j}z}}{2^{j/2}}\cdot\widehat{\psi}\left(\frac z{2^j}\right).
\end{equation}

We begin with the case $k\not= l,$ $k\not= 0,$ $l\not= 0.$

Applying integration by parts,
the assumptions of the theorem, and twice the intermediate value theorem of derivatives $\left|\widehat{\psi}(\frac
z{2^j})\right|=\left|\widehat{\psi}'(\tilde z)\frac
z{2^j}\right|\le c_{\psi''}\frac {|z|^2}{2^{2j}}$ yields the following
 $$|\mathbf E\eta_{jk}\overline{\eta_{jl}}|=\left|\frac{1}{2 \pi}\int\limits_{\mathbb R}\widehat{R}(z)
 \frac{e^{i\frac {k-l}{2^j}z}}{2^{j}}\left|\widehat{\psi}\left(\frac z{2^j}\right)\right|^2\,dz\right|=\left|\frac{1}{2i\pi (k-l)}\right.
\left[\left.\widehat{R}(z)
 e^{i\frac {k-l}{2^j}z}\left|\widehat{\psi}\left(\frac z{2^j}\right)\right|^2\right|_{z=-\infty}^{+\infty}\right.$$
 $$-\int\limits_{\mathbb R}\left(\widehat{R}'(z)\left|\widehat{\psi}\left(\frac z{2^j}\right)\right|^2
 \left.\left.+\widehat{R}(z)\frac{2}{2^j}\,\Re\left(\overline{\widehat{\psi}\left(\frac z{2^j}\right)}
\widehat{\psi}'\left(\frac z{2^j}\right)\right)\right)e^{i\frac
{k-l}{2^j}z}\,dz\right]\right|$$
$$\le\frac{1}{2 \pi|k-l|}\int\limits_{\mathbb R}\left(|\widehat{R}'(z)|\left(c_{\psi''}\frac {|z|^2}{2^{2j}}\right)^2
+2|\widehat{R}(z)|c_{\psi''}\frac {|z|^2}{2^{2j}}\,\frac {|z|}{2^{2j}}
c_{\psi''}\right)\,dz=\frac{A^{\psi}}{2^{4j}|k-l|}\,,
$$
where $$A^{\psi}:=\frac{c_{\psi''}^2}{2 \pi}\int\limits_{\mathbb
R}\left(|\widehat{R}'(z)||z|^4
+2|\widehat{R}(z)||z|^3\right)\,dz<\infty\,,$$
because of conditions \ref{con5} and \ref{con6}.

We use (\ref{104.11}) to estimate the term $|{\psi}_{jl}(t)-{\psi}_{jl}(s)|\,.$ Then
$${\psi}_{jl}(t)=\int\limits_{\mathbb R}\frac {e^{itz}e^{-i\frac l{2^j}z}}{\pi2^{j/2+1}}\widehat{\psi}\left(\frac z{2^j}\right)\,dz=\int\limits_{\mathbb R}\frac{e^{it\left(z+\frac{2^j}l\pi\right)}e^{-i\left(\frac l{2^j}z+\pi\right)}}{\pi 2^{j/2+1}}\widehat{\psi}\left(\frac z{2^j}+\frac {\pi}l\right)\,dz\,$$
$$=\frac{1}{2^{j/2+2}\pi}\int\limits_{\mathbb R}e^{-i\frac l{2^j}z}\left(e^{itz}\widehat{\psi}\left(\frac z{2^j}\right)-e^{it\left(z+\frac{2^j}l\pi\right)}\widehat{\psi}\left(\frac z{2^j}+\frac {\pi}l\right)\right)\,dz\,.
$$
Therefore
$$|{\psi}_{jl}(t)-{\psi}_{jl}(s)|\le\frac{1}{2^{j/2+2}\pi}\int\limits_{\mathbb
R}\left|\left(e^{itz}-e^{isz}\right)\widehat{\psi}\left(\frac
z{2^j}\right)-\left(e^{it\left(z+\frac{2^j}l\pi\right)}-e^{is\left(z+\frac{2^j}l\pi\right)}\right)\right.$$
$$\times\left.\widehat{\psi}\left(\frac
z{2^j}+\frac {\pi}l\right)\right|\,dz\le\frac{2^{j/2-2}}{\pi}\left(\int\limits_{\mathbb
R}\left|e^{it2^ju}-e^{is2^ju}-e^{it2^j\left(u+\frac{\pi}l\right)}+e^{is2^j\left(u+\frac{\pi}l\right)}
\right|\left|\widehat{\psi}(u)\right|du\right.$$
\begin{equation}\label{I2}+\left.\int\limits_{\mathbb
R}\left|e^{it2^j\left(u+\frac{\pi}l\right)}-e^{is2^j\left(u+\frac{\pi}l\right)}
\right|\left|\widehat{\psi}(u)-\widehat{\psi}\left(u+\frac
{\pi}l\right)\right|\,du\right)=:\frac{2^{j/2-2}}{\pi}\,(I_1+I_2)\,.\end{equation}
By the inequality (59) given in \cite{kozroz}
\begin{equation}\label{105}|e^{itz}-e^{isz}|=2\left|\sin\left(\frac{z(t-s)}2\right)\right|
\le2\left(\frac{\ln\left(e^{\alpha}+ \frac{|z|}2\right)}{\ln\left(e^{\alpha}+
\frac1{|t-s|}\right)}\right)^{\alpha}, \quad \alpha>0\,.
\end{equation}
An application of this inequality to the second integral in (\ref{I2}) results~in
$$I_2=\int\limits_{\mathbb
R}\left|e^{it2^jv}-e^{is2^jv} \right|\left|\widehat{\psi}\left(v-\frac
{\pi}l\right)-\widehat{\psi}(v)\right|\,dv\le\frac2{\left(\ln\left(e^{\alpha}+
\frac1{|t-s|}\right)\right)^{\alpha}}$$
$$\times\int\limits_{\mathbb
R}\left(\ln\left(e^{\alpha}+
\frac{2^j|v|}2\right)\right)^{\alpha}\left|\widehat{\psi}\left(v-\frac
{\pi}l\right)-\widehat{\psi}(v)\right|^{\beta}\left|\widehat{\psi}\left(v-\frac
{\pi}l\right)-\widehat{\psi}(v)\right|^{1-\beta}\,dv,$$
where $\beta:=1-\gamma\in (1/2,1).$

In the following derivations, we will use conditions \ref{con1}, \ref{con4}
and the estimates
$$\left|\widehat{\psi}\left(v-\frac{\pi}l\right)-\widehat{\psi}(v)\right|^{\beta}\le c_{\psi'}^\beta\left(\frac{\pi}l\right)^\beta,$$
$$\left|\widehat{\psi}\left(v-\frac
{\pi}l\right)-\widehat{\psi}(v)\right|^{1-\beta}\le 2^{1-\beta}\left(\left|\widehat{\psi}\left(v-\frac
{\pi}l\right)\right|^{1-\beta}+\left|\widehat{\psi}(v)\right|^{1-\beta}\right)\,,$$
\begin{equation}\label{5}\int\limits_{\mathbb
R}\left(\ln\left(e^{\alpha}+
\frac{2^j|v|}2\right)\right)^{\alpha}\left|\widehat{\psi}\left(v\right)\right|^{1-\beta}\,dv\le\int\limits_{\mathbb
R}\left(\ln\left[5^{j+1}\left(\frac{e^{\alpha}}{5^{j+1}}+ \frac{2^{j-1}}{5^{j+1}}\left|v\right|\right)\right]\right)^{\alpha}
\end{equation}
$$\times\left| \widehat{\psi}(v)\right|^{1-\beta}\,dv\le 2^\alpha\left((\ln5)^\alpha (j+1)^{\alpha}c_0+c_1\right)<\infty\,,$$

$$\int\limits_{\mathbb
R}\left(\ln\left(e^{\alpha}+
\frac{2^j|v|}2\right)\right)^{\alpha}\left|\widehat{\psi}\left(v-\frac
{\pi}l\right)\right|^{1-\beta}\,dv\le\int\limits_{\mathbb
R}\left(\ln\left[5^{j+1}\left(\frac{e^{\alpha}}{5^{j+1}}+ \frac{2^{j-1}}{5^{j+1}}\left(\left|v\right|+\frac{\pi}{|l|}\right)\right)\right]\right)^{\alpha}$$
$$\times \left|\widehat{\psi}(v)\right|^{1-\beta}\,dv
\le 2^\alpha\left((\ln5)^\alpha (j+1)^{\alpha}c_0+c_1\right)<\infty\,,$$
where
$$c_0:=\int\limits_{\mathbb
R}\left|\widehat{\psi}(v)\right|^{1-\beta}\,dv<\infty\,,\quad c_1:=\int\limits_{\mathbb
R}\left(\ln(1+|v|)\right)^{\alpha}\left|\widehat{\psi}(v)\right|^{1-\beta}\,dv<\infty\,.$$

The integral $c_0$ is finite because of the boundedness of $\widehat{\psi}(v)$ and condition~\ref{con4}.

Using these facts, we get
  \begin{equation}\label{205}I_2\le\frac{2^{3+\alpha-\beta}\pi^\beta c_{\psi'}^\beta}{|l|^{\beta}\left(\ln\left(e^{\alpha}+
\frac1{|t-s|}\right)\right)^{\alpha}}\left((\ln5)^\alpha (j+1)^{\alpha}c_0+c_1\right)\,.\end{equation}

Similarly, we can  estimate the first integral. It is easy to see that

$$|\Delta|:=\left|e^{it2^ju}-e^{is2^ju}-e^{it2^j(u+\frac{\pi}l)}+e^{is2^j(u+
\frac{\pi}l)}\right|\le\left|e^{it2^ju}-e^{is2^ju}\right|\left|1-e^{it2^j\frac{\pi}l}\right|$$
$$+\left|e^{is2^j\frac {\pi}l}-e^{it2^j\frac{\pi}l}\right|\le2\left|\sin\left(\frac{2^ju(t-s)}2\right)\right|\cdot\left|\sin\left(\frac{2^j\pi t}{2l}\right)\right|+2\left|\sin\left(\frac{2^j\pi(t-s)}{2l}\right)\right|\,.$$
Note that, by (\ref{105}) and (\ref{5}):
\begin{equation}\label{+}\left|\sin\left(\frac{2^ju(t-s)}2\right)\right|\le\left(\frac{\ln(e^{\alpha}+ \frac{2^j|u|}2)}{\ln\left(e^{\alpha}+
\frac1{|t-s|}\right)}\right)^{\alpha}
\le 2^\alpha\frac{(j+1)^\alpha (\ln 5)^\alpha+\left(\ln\left(1+ |u|\right)\right)^\alpha}{\left(\ln\left(e^{\alpha}+
\frac1{|t-s|}\right)\right)^{\alpha}}\,.\end{equation}
Due to \cite[Lemma 4.2]{koztur}

\begin{equation}\label{*}\left|\sin\left(\frac{2^j\pi(t-s)}{2l}\right)\right|\le \left|\frac{2^j\pi(t-s)}{2l}\right|\le \frac{2^jc_{\alpha}}{|l|\left(\ln\left(e^{\alpha}+ \frac1{|t-s|}\right)\right)^{\alpha}},\quad  \alpha>0\,,\end{equation}
where $c_{\alpha}$ depends only on $T $ and $\alpha\,.$

Applying inequalities (\ref{+}), (\ref{*}) and $\left|\sin\left(\frac{2^j\pi t}{2l}\right)\right|\le \frac{2^{j-1}\pi T}{|l|}$
 we get
$$|\Delta|\le\frac{2^{j+1}\left(\pi T 2^{\alpha-1}\left[(j+1)^\alpha (\ln 5)^\alpha+\left(\ln\left(1+ |u|\right)\right)^\alpha\right]+c_{\alpha}\right)}{|l|\left(\ln\left(e^{\alpha}+ \frac1{|t-s|}\right)\right)^{\alpha}}
\,.$$

Using above inequalities the first integral can be estimated as follows:
$$I_1 \le\frac{2^{j+1}}{|l|\left(\ln\left(e^{\alpha}+
\frac1{|t-s|}\right)\right)^{\alpha}}\int\limits_{\mathbb R}\left|\widehat{\psi}(u)\right|
\left(\pi T 2^{\alpha-1}\left((j+1)^\alpha (\ln 5)^\alpha+c_{\alpha}\right.\right. $$
\begin{equation}\label{206}
\left.\left.+\left(\ln\left(1+ |u|\right)\right)^\alpha\right)\right)\,du=\frac{2^{j+1}\left(\pi T 2^{\alpha-1}\left((j+1)^\alpha (\ln 5)^\alpha c_2+c_3\right)+c_{\alpha}c_2\right)}{|l|\left(\ln\left(e^{\alpha}+
\frac1{|t-s|}\right)\right)^{\alpha}}\,,\end{equation}
where
$$c_2:=\int\limits_{\mathbb
R}\left|\widehat{\psi}(v)\right|\,dv<\infty\,,\quad c_3:=\int\limits_{\mathbb
R}\left(\ln(1+|v|)\right)^{\alpha}\left|\widehat{\psi}(v)\right|\,dv<\infty\,.$$

The integrals $c_2$ and $c_3$ are finite because $\widehat{\psi}(v)$ is bounded:

\begin{equation}\label{c2}c_2\le\sup\limits_{u\in\mathbb R}\left|\widehat{\psi}(u)\right|\int\limits_{\mathbb
R}\frac{\left|\widehat{\psi}(v)\right|^{1-\beta}}{\left(\sup_{u\in\mathbb R}\left|\widehat{\psi}(u)\right|\right)^{1-\beta}}\,dv=
\left(\sup\limits_{u\in\mathbb R}\left|\widehat{\psi}(u)\right|\right)^{\beta}c_0<\infty\,,
\end{equation}
$$c_3\le \left(\sup\limits_{u\in\mathbb R}\left|\widehat{\psi}(u)\right|\right)^{\beta}c_1<\infty\,.$$

Using (\ref{205}) and (\ref{206}), we obtain:
$$|{\psi}_{jk}(t)-{\psi}_{jk}(s)|\cdot|{\psi}_{jl}(t)-{\psi}_{jl}(s)|\le
\frac{2^{j-4}}{|k|^{\beta}|l|^{\beta}}\left(\frac{2^{3+\alpha-\beta}\pi^\beta c_{\psi'}^\beta\left((\ln5)^\alpha (j+1)^{\alpha}c_0+c_1\right)}{\pi \left(\ln\left(e^{\alpha}+
\frac1{|t-s|}\right)\right)^{\alpha}}\right.$$
\begin{equation}\label{ff}+\left.\frac{2^{j+1}\left(\pi T 2^{\alpha-1}\left((j+1)^\alpha  (\ln 5)^\alpha c_2+c_3\right)+c_{\alpha}c_2\right)}{\pi \left(\ln\left(e^{\alpha}+
\frac1{|t-s|}\right)\right)^{\alpha}} \right)^2\le
\frac{(j+1)^{2\alpha}2^{3j-2}K^2}{|k|^{\beta}|l|^{\beta}\left(\ln\left(e^{\alpha}+
\frac1{|t-s|}\right)\right)^{2\alpha}}\,,
\end{equation}
where
$$K:= \pi^{-1}\left(2^{3+\alpha-\beta}\pi^\beta c_{\psi'}^\beta\left((\ln5)^\alpha c_0+c_1\right)+\pi T 2^{\alpha-1}\left((\ln 5)^\alpha c_2+c_3\right)+c_{\alpha}c_2\right)\,.$$

Thus
\begin{equation}\label{knel}\sum\limits_{\scriptsize\begin{array}{c}
                     |k|\le k_j,|l|\le k_j\\
                      k\ne l,kl\ne0
                   \end{array}}
                   \hspace{-5mm}|\mathbf E\eta_{jk}\overline{\eta_{jl}}||\psi_{jk}(t)-\psi_{jk}(s)|
|\psi_{jl}(t)-\psi_{jl}(s)|\le \frac{(j+1)^{2\alpha}A^{\psi}QK^2}{2^{j}\left(\ln\left(e^{\alpha}+
\frac1{|t-s|}\right)\right)^{2\alpha}} \,,
\end{equation}
where

$$Q:=\sum\limits_{\scriptsize\begin{array}{c}
                     |k|\le k_j,|l|\le k_j\\
                      k\ne l,kl\ne0
                   \end{array}}
\frac{1}{4|k-l||k|^{\beta}|l|^{\beta}}\le\sum\limits_{0< k\le k_j} \sum\limits_{0<l\le k_j}\frac {1}{4\sqrt {kl}k^{\beta}l^{\beta}} +\sum\limits_{0< m< k_j} \sum\limits_{0<l\le
k_j-m}\frac {1}{ml^{\beta}(m+l)^{\beta}}\,.$$

Using the inequality
$x+y>c_{\delta}x^{\delta}y^{1-\delta}, {\ }x,y>0, {\
}\delta\in(0,1),$ $c_{\delta}:={\delta^{-\delta}(1-\delta)^{\delta-1}}$, we obtain
$$Q\le\sum\limits_{0< k\le k_j} \sum\limits_{0<l\le k_j}\frac {1}{4k^{\frac12+\beta}l^{\frac12+\beta}}
 +c^\beta_{\delta}\sum\limits_{0< m< k_j}
\sum\limits_{0<l\le k_j-m}\frac {1}{ml^{\beta}m^{\delta\beta}l^{(1-\delta)\beta}}$$
$$<\left(\sum\limits_{k=1}^\infty \frac {1}{2k^{\frac12+\beta}}\right)^2
+c^\beta_{\delta}\sum\limits_{m=1}^\infty\frac {1}{m^{1+\delta\beta}} \sum\limits_{l=1}^\infty\frac
{1}{l^{(2-\delta)\beta}}<\infty\,.$$
 The statement becomes apparent if $\delta\in (0,2-1/\beta)$ is chosen.

Similarly, we can  exploit the case $l\ne0, k=0:$
$$|\mathbf E\eta_{j0}\overline{\eta_{jl}}|=\frac{1}{2 \pi }\left|\int\limits_{\mathbb R}\widehat{R}(z)
\cdot\frac{e^{-i\frac
l{2^j}z}}{2^{j}}\left|\widehat{\psi}\left(\frac z{2^j}\right)\right|^2\,dz\right|=\left|-\frac{1}{2i\pi l}\right.\left[\left.\widehat R(z)
e^{-i\frac l{2^j}z}\left|\widehat{\psi}\left(\frac
z{2^j}\right)\right|^2\right|_{z=-\infty}^{+\infty}\right.$$
$$-\int\limits_{\mathbb R}\left(\widehat{R}'(z)\left|\widehat{\psi}\left(\frac z{2^j}\right)\right|^2
\left.\left.+\widehat R(z)\frac{2}{2^j}\,\Re\left(\widehat{\psi}\left(\frac
z{2^j}\right)\widehat{\psi}'\left(\frac z{2^j}\right)\right)\right) e^{-i\frac
l{2^j}z}\,dz\right]\right|$$
\begin{equation}\label{eta}\le\frac{1}{2 \pi|l|}\int\limits_{\mathbb R}\left(|\widehat{R}'(z)|\left(c_{\psi''}\frac {|z|^2}{2^{2j}}\right)^2
+2|\widehat{R}(z)|c_{\psi''}\frac {|z|^2}{2^{2j}}\frac {|z|}{2^{j}}
c_{\psi''}\right)\,dz=\frac{A^{\psi}}{2^{4j}|l|}\,.
\end{equation}

By (\ref{105}) and (\ref{5})
$$|{\psi}_{j0}(t)-{\psi}_{j0}(s)|=\frac{2^{j/2}}{2
\pi}\left|\int\limits_{\mathbb R}(e^{it2^{j}u}-e^{is2^{j}u})
\widehat{\psi}(u)\,du\right|\le\frac{2^{j/2}}{
\pi}\int\limits_{\mathbb
R}\left(\frac{\ln\left(e^{\alpha}+\frac{2^j|u|}2\right)}{\ln\left(e^{\alpha}+
\frac1{|t-s|}\right)}\right)^{\alpha}|\widehat{\psi}(u)|\,du$$
\begin{equation}\label{0}\le \int\limits_{\mathbb
R}\frac{2^{j/2}\left((j+1)\ln 5+\ln(1+|u|)\right)^{\alpha}}{
\pi\left(\ln\left(e^{\alpha}+
\frac1{|t-s|}\right)\right)^{\alpha}}|\widehat{\psi}(u)|\,du
\le\frac{2^{j/2+\alpha}\left((j+1)^\alpha\left(\ln 5\right)^\alpha c_2+c_3\right)}{\pi\left(\ln\left(e^{\alpha}+ \frac1{|t-s|}\right)\right)^{\alpha}}\,.\end{equation}

Then
\begin{equation}\label{105.1}\sum\limits_{0<|l|\le k_j}|\mathbf E\eta_{j0}\overline{\eta_{jl}}||\psi_{j0}(t)-\psi_{j0}(s)||\psi_{jl}(t)-\psi_{jl}(s)|
\le\frac{(j+1)^{2\alpha}\cdot q}{2^{2j}\left(\ln\left(e^{\alpha}+
\frac1{|t-s|}\right)\right)^{2\alpha}}\,,
\end{equation}
where
$$q:=\frac{2^{\alpha}A^{\psi}K((\ln 5)^\alpha c_2+c_3)}{\pi}\cdot\sum\limits_{l=1}^\infty
\frac{1}{|l|^{1+\beta}}<\infty\,.$$

For $l=k, k\not=0$ we get the estimates

$$|{\psi}_{jk}(t)-{\psi}_{jk}(s)|^2\le
\frac{(j+1)^{2\alpha}2^{3j-2}K^2}{|k|^{2\beta}\left(\ln\left(e^{\alpha}+
\frac1{|t-s|}\right)\right)^{2\alpha}}\,,$$

\begin{equation}\label{kl}\mathbf E|\eta_{jk}|^2\le\frac{1}{2^{j+1} \pi}
\int\limits_{\mathbb R}|\widehat{R}(z)|\left(c_{\psi''}\frac {|z|^2}{2^{2j}}\right)^2\,dz=\frac{A^{\psi}_1}{2^{5j}}\,,\end{equation}
where $$A^{\psi}_1:=\frac{c_{\psi''}^2}{2 \pi}\int\limits_{\mathbb
R}|\widehat{R}(z)||z|^4\,dz<\infty\,.$$

Hence
\begin{equation}\label{kel}\sum\limits_{k=1}^{k_j} \mathbf E|\eta_{jk}|^2|{\psi}_{jk}(t)-{\psi}_{jk}(s)|^2
\le \frac{(j+1)^{2\alpha}\cdot q_1}{2^{2j+1}\left(\ln\left(e^{\alpha}+
\frac1{|t-s|}\right)\right)^{2\alpha}}\,,
\end{equation}
where $$q_1:= \frac{A^{\psi}_1 K^2}{2}\cdot\sum\limits_{k=1}^{\infty}\frac{1}{|k|^{2\beta}}<\infty.$$

Finally, for $k=l=0,$ applying (\ref{kl}) and (\ref{0})  we get

\begin{equation}\label{kl0}
\mathbf E|\eta_{j0}|^2\cdot|{\psi}_{j0}(t)-{\psi}_{j0}(s)|^2\le \frac{(j+1)^{2\alpha}\cdot q_2}{2^{4j}\left(\ln\left(e^{\alpha}+ \frac1{|t-s|}\right)\right)^{2\alpha}}\,, \end{equation}
where
$$q_2:=\frac{2^{2\alpha}A^{\psi}_1}{\pi^2} \left(\left(\ln 5\right)^\alpha c_2+c_3\right)^2\,.$$

Using (\ref{knel}), (\ref{105.1}),  (\ref{kel}) and (\ref{kl0}) we find that
$$\sum\limits_{j=0}^{n-1}\sqrt{S_j}\le\sum\limits_{j=0}^{n-1}\left(
\frac{(j+1)^{2\alpha}\cdot q_2}{2^{4j}\left(\ln\left(e^{\alpha}+ \frac1{|t-s|}\right)\right)^{2\alpha}}+
\frac{(j+1)^{2\alpha}\cdot q_1}{2^{2j}\left(\ln\left(e^{\alpha}+
\frac1{|t-s|}\right)\right)^{2\alpha}}+\right.$$
\begin{equation}\label{1066}\left.
\frac{(j+1)^{2\alpha}A^{\psi}QK^2}
{2^{j}\left(\ln\left(e^{\alpha}+ \frac1{|t-s|}\right)\right)^{2\alpha}}+
\frac{2(j+1)^{2\alpha}\cdot q}{2^{2j}\left(\ln\left(e^{\alpha}+
\frac1{|t-s|}\right)\right)^{2\alpha}}\right)^{1/2}
\le \frac{B_1}{\left(\ln\left(e^{\alpha}+
\frac1{|t-s|}\right)\right)^{\alpha}}\,,
\end{equation}
where
$$B_1:=\left(
q_2+q_1+A^{\psi}QK^2+2q\right)^{1/2}\cdot\sum\limits_{j=0}^{\infty}\frac{(j+1)^{\alpha}}{2^{j/2}}<\infty\,.$$

The results of identical analysis for $S$ are given below.

First, we evaluate  $|\mathbf E\xi_{0k}\overline{\xi_{0l}}|$ in the case $k\not= l,$ $k\ne0,$ $l\ne0\,:$

$$|\mathbf E\xi_{0k}\overline{\xi_{0l}}|=\left|\frac{1}{2
\pi}\int\limits_{\mathbb R}\widehat
R(z)e^{i(k-l)z}\left|\widehat{\phi}( z)\right|^2 \,dz\right|\le \frac{ A^{\phi}}{|k-l|}\,,
$$
where
$$A^{\phi}=\frac{1}{2 \pi}\,\left(c_{\phi}^2\int\limits_{\mathbb
R}\left|\widehat R'(z)\right|\,dz
+2c_{\phi}c_{\phi'}\int\limits_{\mathbb R}\left|\widehat
R(z)\right|\,dz\right)<\infty\,$$
because of assumptions \ref{con5} and \ref{con6}.

Similarly to (\ref{ff}) we  derive
$$|\phi_{0k}(t)-\phi_{0k}(s)|\cdot|\phi_{0l}(t)-\phi_{0l}(s)|
\le \frac {(K^{\phi})^2}{4|l|^{\beta}|k|^{\beta}\left(\ln\left(e^{\alpha}+
\frac1{|t-s|}\right)\right)^{2\alpha}}\,,
$$ where
$$K^{\phi}=\pi^{-1}\left(2^{3+\alpha-\beta}\pi^\beta c_{\phi'}^\beta\left((\ln5)^\alpha c_{\phi 0}+c_{\phi 1}\right)+\pi T 2^{\alpha-1}\left((\ln 5)^\alpha c_{\phi 2}+c_{\phi 3}\right)+c_{\alpha}c_{\phi 2}\right)\,,$$
$$c_{\phi 0}:=\int\limits_{\mathbb
R}\left|\widehat{\phi}(v)\right|^{1-\beta}\,dv<\infty\,,\quad c_{\phi 1}:=\int\limits_{\mathbb
R}\left(\ln(1+|v|)\right)^{\alpha}\left|\widehat{\phi }(v)\right|^{1-\beta}\,dv<\infty\,,$$

$$c_{\phi  2}:=\int\limits_{\mathbb
R}\left|\widehat{\phi }(v)\right|\,dv<\infty\,,\quad c_{\phi  3}:=\int\limits_{\mathbb
R}\left(\ln(1+|v|)\right)^{\alpha}\left|\widehat{\phi }(v)\right|\,dv<\infty\,.$$

Therefore
\begin{equation}\label{1070}\sum\limits_{\scriptsize\begin{array}{c}
                     |k|\le k_0,|l|\le k_0\\
                      k\ne l,kl\ne0
                   \end{array}}|
\mathbf E\xi_{0k}\overline{\xi_{0l}}||\phi_{0k}(t)-\phi_{0k}(s)|
|\phi_{0l}(t)-\phi_{0l}(s)|
\le
\frac {A^{\phi}(K^{\phi})^2 Q}{\left(\ln\left(e^{\alpha}+
\frac1{|t-s|}\right)\right)^{2\alpha}}\,.
\end{equation}
Similarly to (\ref{eta}) for the case $k=0,$ $l\ne0$ we obtain:
$$|\mathbf E\xi_{00}\overline{\xi_{0l}}|\le\frac{1}{2 \pi|l|}\int\limits_{\mathbb R}\left(\left|\widehat
R'(z)\right|\,|\widehat{\phi}(z)|^2+2\left|\widehat R(z)\right|\,|\overline{\widehat{\phi}(z)}
\widehat{\phi}'(z)|\right)\,dz\le\frac{A^{\phi}}{|l|}
\,.$$

Similarly to (\ref{0}), by (\ref{105}) we obtain
\begin{equation}\label{f00}|{\phi}_{00}(t)-{\phi}_{00}(s)|
\le\frac1{\pi}\int\limits_{\mathbb R}
\left(\frac{\ln\left(e^{\alpha}+ \frac{|z|}2\right)}{\ln\left(e^{\alpha}+
\frac1{|t-s|}\right)}\right)^{\alpha}|\widehat{\phi}(z)|\,dz\le\frac{2^\alpha ((\ln 5)^{\alpha}c_{\phi 2}+c_{\phi 3})}{\pi\left(\ln\left(e^{\alpha}+
\frac1{|t-s|}\right)\right)^{\alpha}}\,.
\end{equation}

Therefore
\begin{equation}\label{1072}\sum\limits_{|l|\le k_0, l\ne0}|\mathbf E\xi_{00}\overline{\xi_{0l}}|\cdot|\phi_{00}(t)-\phi_{00}(s)|\cdot
|\phi_{0l}(t)-\phi_{0l}(s)|\le\frac{q_\phi}{\left(\ln\left(e^{\alpha}+
\frac1{|t-s|}\right)\right)^{2\alpha}}\,,
\end{equation}
where
$$q_\phi:=2^\alpha \pi^{-1}A^{\phi} K^\phi \left((\ln 5)^{\alpha}c_{\phi 2}+c_{\phi 3}\right)\cdot\sum_{l=1}^\infty \frac{1}{|l|^{1+\beta}}<\infty\,.$$

For $l=k,$ $k\not= 0$ we obtain
$$|\phi_{0k}(t)-\phi_{0k}(s)|^2\le \frac {(K^{\phi})^2}{4\,|k|^{2\beta}\left(\ln\left(e^{\alpha}+
\frac1{|t-s|}\right)\right)^{2\alpha}}\,,$$

\begin{equation}\label{Af}\mathbf E|\xi_{0k}|^2=\frac{1}{2 \pi}\int\limits_{\mathbb
R}\widehat R(z)|\widehat{\phi}( z)|^2dz\,\le \frac{c_\phi^2}{2\pi}\int\limits_{\mathbb
R}|\widehat R(z)|\,dz=:A_1^\phi<\infty\,.
\end{equation}

Hence

\begin{equation}\label{kel1}\sum\limits_{k=1}^{k_0} \mathbf E|\xi_{0k}|^2|{\phi}_{0k}(t)-{\phi}_{0k}(s)|^2
< \frac{q_{\phi 1}}{2\left(\ln\left(e^{\alpha}+
\frac1{|t-s|}\right)\right)^{2\alpha}}\,,
\end{equation}
where $$q_{\phi 1}:= \frac{ A^{\phi}_1 (K^\phi)^2}{2}\cdot\sum\limits_{k=1}^{\infty}\frac{1}{|k|^{2\beta}}<\infty.$$

Finally, for $k=l=0,$ by (\ref{Af}) and (\ref{f00}) we get
\begin{equation}\label{kel0}\mathbf E|\xi_{00}|^2\cdot |\phi_{00}(t)-\phi_{00}(s)|^2\le \frac{q_{\phi 2}}{\left(\ln\left(e^{\alpha}+
\frac1{|t-s|}\right)\right)^{2\alpha}}\,,
\end{equation}
 where $$q_{\phi 2}:= \frac{2^{2\alpha} A^{\phi}_1}{\pi^2}((\ln 5)^\alpha c_{\phi 2} +c_{\phi 3})^2\,.$$

Using (\ref{1070}), (\ref{1072}),  (\ref{kel1}), and (\ref{kel0}) we have
\begin{equation}\label{1074}\sqrt{S}\le\frac {B_2}{\left(\ln\left(e^{\alpha}+
\frac1{|t-s|}\right)\right)^{\alpha}},
\end{equation}
where
$$B_2:=\left(q_{\phi 1}+q_{\phi 2}+A^{\phi}(K^{\phi})^2 Q+2 q_\phi\right)^{1/2}\,.$$

By combining (\ref{1066}) and (\ref{1074}) we obtain that (\ref{1075}) is valid for $B=B_1+B_2\,.$

Using Theorem~\ref{8103} we get $\mathbf{X}_{n,\mathbf{k}_n}(t)\to
\mathbf{X}(t)$ in mean square when $n\to\infty$ and $k_j\to\infty$ for all $j\in N_0\,.$
Now, combining  (\ref{1075}), Theorem~\ref{213} and
Remark~\ref{772}, we conclude that $\mathbf{X}_{n,\mathbf{k}_n}(t)\to \mathbf{X}(t)$ uniformly in probability on each
interval $[0,T]$ when $n\to\infty$ and $k_j\to\infty$ for all $j\in \mathbb N_0\,.$

\section{Conclusions}
We have analyzed uniform convergence of wavelet
expansions of stationary Gaussian random processes. The most general form of the expansions was studied.
The results are obtained under simple conditions which can be easily verified. The conditions  are weaker than those in the former literature. The main theorem is the first result on stochastic uniform convergence of general finite wavelet expansions in the open literature.

The wavelet expansions are more general than the Fourier-wavelets decompositions studied in \cite{kur}. The integrals $\alpha_{0k}$ and $\beta_{jk}$ may be dependent, which can be used to derive similar uniform results for nonstationary random processes (see, for example, \cite{mey}).

It would be of interest to adopt the results to various wavelet bases, which are important for applications.

By Theorem~\ref{8103}, it is possible to obtain wavelet expansions (\ref{Xn}) for numerous wide classes of random processes (for example, the fractional Brownian motion) in $L_2.$  The technique, which was developed in Theorem~\ref{main}, can be used to prove uniform convergence. However, the generalizations to other random processes are not always straightforward, because conditions \ref{con5} and \ref{con6} may not be satisfied. The essence of the problem is a modification of the proof by obtaining new upper bounds on $|\mathbf E\eta_{jk}\overline{\eta_{jl}}|$ and $|\mathbf E\xi_{0k}\overline{\xi_{0l}}|.$ For some classes of stochastic processes this non-trivial problem requires a full-length paper itself. The authors are preparing such results for the fractional Brownian motion now.

Some new results on the rate of convergence of the wavelet expansions in the space $C([0,T])$ were obtained too. These results will be published in a separate paper.

\section{Acknowledgements}
This work was partly supported by
La Trobe University Research Grant--501821 "Sampling, wavelets and optimal stochastic modelling".

\end{document}